\documentclass[]{article}
\begin{document}
\newtheorem{proposition}{Proposition}[section]
\newtheorem{definition}{Definition}[section]
\newtheorem{lemma}{Lemma}[section]

\title{\bf Hopf-like Algebras\\ and  \\Extended P-B-W Theorems}
\author{\small{Keqin Liu}\\\small{Department of Mathematics, The University of British Columbia}\\\small{Vancouver, BC,
Canada, V6T 1Z2}}
\date{\small{November, 2010}}
\maketitle

{\bf Abstract} {\small Based on invariant algebras, we introduce representations$^{6-th}$ of Lie algebras and representations$^{\langle 4-th\rangle }$ of Leibniz algebras, give the  extended P-B-W Theorems in the context of the new representations of Lie algebras and Leibniz algebras, and generalize the Hopf-algebra structure on the enveloping algebras of Lie Algebras.}

\bigskip
Throughout, an associative algebra always means an associative algebra having an identity, a homomorphism from an associative algebra to an associative algebra always preserves the identity, and all vector spaces are vector spaces over a field $\mathbf{k}$.

\bigskip
Let $V$ be a vector space, and let $End(V)$ be the associative algebra of all linear transformations from $V$ to $V$. It is well-known that $End(V)$ is a Lie algebra with respect to the following bracket $[\, ,\,]$
\begin{equation}\label{eq1}
[f, g]:=fg-gf \quad\mbox{for $f$, $g\in End(V)$}.
\end{equation}
This Lie algebra is denoted by $g\ell (V)$. In this paper, the bracket defined by (\ref{eq1}) will be called the {\bf ordinary bracket}, and a homomorphism from a Lie algebra $\mathcal{L}$ to the Lie algebra  $g\ell (V)$ will be called an {\bf ordinary representation} of the Lie algebra $\mathcal{L}$ on $V$.

\medskip
This paper is concerned with the new ways of representing Lie algebras and Leibniz algebras by using linear transformations. The key idea in this new ways is to replace the ordinary bracket by new brackets and replace the associative algebra $End(V)$ by a subalgebra of $End(V)$. Let $W$ be a non-zero subspace of a vector space $V$, and let $End _{_W}(V)$ be the set   of all  linear transformations from $V$ to $V$ which have $W$ as an invariant subspace, i.e.,
$$
End _{_W}(V):=\{\, f \,|\, \mbox{$f\in End(V)$ and $f(W)\subseteq W$}\,\}.
$$
Clearly, $End _{_W}(V)$ is a subalgebra of the associative algebra $End(V)$. Although  the way of making $End(V)$ into a Leibniz algebra has not been found, there are several ways of making the subalgebra $End _{_W}(V)$ into a Leibniz algebra. Also, there are several new brackets which are different from the ordinary bracket and make $End _{_W}(V)$ into a Lie algebra. In fact, if we fix a scalar $k\in\mathbf{k}$ and a linear transformation $q$ satisfying
\begin{equation}\label{eq2}
q(W)=0\quad\mbox{and}\quad\mbox {$q(v)-v\in W$ for all $v\in V$},
\end{equation}
then $End _{_W}(V)$ becomes a Leibniz algebra  with respect to the following $4$-th angle bracket
$\langle \,, \, \rangle _{4,k}$:
$$
\langle f, g \rangle _{4, k} =fg-gf+gqf-fgq+ kfqg - kqgf\quad\mbox{for $f$, $g\in End _{_W}(V)$},
$$
and $End _{_W}(V)$ becomes a Lie algebra with respect to the following $6$-th square bracket $[\, ,\,]_{6, k}$:
$$
[f, g]_{6, k}:=fg-gf-fgq+gfq+ kfqg - kgqf\quad\mbox{for $f$, $g\in End _{_W}(V)$}.
$$
Replacing $g\ell (V)$ by the Leibniz algebra 
($End _{_W}(V)$, $\langle \,, \, \rangle _{4,k}$), we get the notion of 
representations$^{\langle 4-th\rangle }$ of Leibniz algebras.  Replacing $g\ell (V)$ by the Lie algebra ($End _{_W}(V)$, $[\, ,\,]_{6, k}$), we get  the notion of representations$^{6-th}$ of Lie algebras.

\medskip
The Poincar$\acute{e}$-Birkhoff-Witt Theorem (P-B-W Theorem) and the Hopf algebra structure for the enveloping algebras of Lie algebras are of great importance to the ordinary representations of Lie algebras. To initiate the study of the new representations of Lie algebras and Leibniz algebras, it is natural to study the counterparts of the P-B-W Theorem and the Hopf algebra structure for the enveloping algebras of Lie algebras in the context of the new representations of Lie algebras and Leibniz algebras. The purpose of this paper is to present our results in this study.

\medskip
This paper consists of six sections.
In section 1,  we introduce the notion of invariant algebras, discuss some basic properties of invariant algebras, and define representations$^{6-th}$ of Lie algebras and 
representations$^{4-th}$ of Leibniz algebras. In section 2, we introduce free invariant algebras and construct a basis for a free invariant algebras. In section 3, we introduce 
enveloping$^{6-th}$ algebras of Lie algebras and give the  extended$^{6-th}$ P-B-W theorem. In section 4, we introduce bialgebras with $\sigma$-counit and Hopf-like$^{6-th}$ algebras. The main results of this section is that the enveloping$^{6-th}$ algebra of a Lie algebra is a Hopf-like$^{6-th}$ algebra. In section 5, we introduce enveloping$^{\langle 4-th \rangle}$ algebras of Leibniz algebras and give the  extended$^{\langle 4-th \rangle}$ P-B-W theorem. In section 6, we introduce Hopf-like$^{{\langle 4-th \rangle}_i}$  algebras with $i=1$, $2$ and study  Hopf-like$^{{\langle 4-th \rangle}_i}$  algebra structures on the enveloping$^{\langle 4-th \rangle}$ algebra of a Leibniz algebra.

\bigskip
\section{Invariant Algebras}

We begin this section by introducing the notion of invariant algebras.

\begin{definition}\label{def1.1} Let $A$ be an associative algebra with an idempotent $q$. The set
\begin{equation}\label{eq3}
(A, q):=\{\, x, \,|\, \mbox{$x\in A$ and $qxq=qx$}\,\}
\end{equation}
is called the {\bf (right) invariant algebra} induced by the idempotent $q$. 
\end{definition}

\medskip
The most important example of invariant algebras is the linear invariant algebra over a vector space. Let $V$ be a vector space, let $End(V)$ be the associative algebra of all linear transformations from $V$ to $V$, and let $I$ be the identity linear transformation of $V$. If $W$ is a subspace of $V$, then an linear transformation $q$ satisfying (\ref{eq2})
is an idempotent. An element $q$ of $End(V)$ satisfying (\ref{eq2}) is called a 
{\bf $W$-idempotent}. 

\begin{definition}\label{def1.2} If $W$ is a subspace of a vector space $V$ and $q$ is a $W$-idempotent, then the invariant algebra $(End(V), q)$ induced by $q$ is called the {\bf linear invariant algebra over $V$} induced by the $W$-idempotent $q$.
\end{definition}

Clearly, the linear invariant algebra $(End(V), q)$ over $V$ induced by a $W$-idempotent $q$  consists of all linear transformations of $V$ having $W$ as an invariant subspace, i.e., 
$(End(V), q)=End _{_W}(V)$.

\medskip
\begin{definition}\label{def1.3} Let $(A,q_{_A})$ and $(B,q_{_B})$ be invariant algebras. A linear map $\phi: (A,q_{_A}) \to (B,q_{_B})$ is called an {\bf invariant homomorphism} if
$$
\phi(xy)=\phi(x)\phi(y)\quad\mbox{for $x$, $y\in (A,q_{_A})$}
$$
and
$$ \phi(1_{_A})=1_{_B}, \qquad \phi(q_{_A})=q_{_B},$$
where $1_{_A}$ and $1_{_B}$ are the identities of $A$ and $B$, respectively. A bijective invariant homomorphism is called an {\bf invariant isomorphism}.
\end{definition}

\medskip
The next proposition shows that any invariant algebra can be imbedded in a linear invariant algebra.

\begin{proposition}\label{pr1.1} If $(A,q)$ is an invariant algebra, then the map
$$\phi: a\mapsto a_{_L} \quad\mbox{for $a\in (A,q)$}$$
is an injective invariant homomorphism from $(A,q)$ to the linear invariant algebra 
$(End(A,q), q_{_L})$ over $(A,q)$ induced by the $(A, q)^{ann}$-idempotent $q_{_L}$, where $a_{_L}$ is the left multiplication defined by
$ a_{_L}(x):=ax \quad\mbox{for $x\in (A.q)$}$
and 
$$(A, q)^{ann}:= \{\, qx-x \,|\, x\in (A, q)\,\}.$$
\end{proposition}

\medskip
There are six ways of making an invariant algebra into a  Lie algebra. One of the six ways is given in the following

\medskip
\begin{proposition}\label{pr1.2} If $q$ is an idempotent of an associative algebras $A$, then $(A,q)$ becomes a Lie algebra under the following {\bf $6$-th square bracket}
\begin{equation}\label{eq4}
[x, y ]_{6, k} := xy-yx-xyq+yxq+kxqy-kyqx,
\end{equation}
where $x$, $y\in (A,q)$, and $k$ is a fixed scalar in the field $\mathbf{k}$. This Lie algebra is denoted by $Lie\Big((A,q), [\, , \, ]_{6, k}\Big)$.
\end{proposition}

\medskip
Following \cite{Loday}, a vector space $\mathcal{L}$ is called a {\bf (right) Leibniz algebra} if there exists a binary operation $\langle \, , \, \rangle$: $\mathcal{L}\times \mathcal{L}\to \mathcal{L}$  such that the {\bf (right) Leibniz identity} holds:
$\langle\langle x, y\rangle , z\rangle=\langle x, \langle y, z\rangle\rangle+ 
\langle\langle x,  z\rangle , y\rangle$ for $x,y,z\in \mathcal{L}$. 

\medskip
There are four ways of making an invariant algebra into a Leibniz algebra. One of the four ways is given in the following

\medskip
\begin{proposition}\label{pr1.3} If $q$ is an idempotent of an associative algebras $A$, then $(A,q)$ becomes a Leibniz algebra under the following {\bf $4$-th angle bracket}
\begin{equation}\label{eq5}
\langle x, y \rangle _{4} =xy-yx+yqx-xyq+ kxqy - kqyx,
\end{equation}
where $x$, $y\in (A,q)$, and $k$ is a fixed scalar in the field $\mathbf{k}$. This Leibniz algebra is denoted by $Leib\Big((A, q), \langle\, , \, \rangle_{4, k}\Big)$.
\end{proposition}

\medskip
Let $W$ be a subspace of a vector space $V$, and let $q$ be a $W$-idempotent. The Lie algebra $Lie\Big((End(V), q), [\, , \, ]_{6, k}\Big)$ is denoted by 
$g\ell _{q, W}^{[ \,, \, ] _{6,k}}(V)$ and is called the {\bf $6$-th general linear Lie algebra} induced by $(q, W)$. The Leibniz algebra 
$Lie\Big((End(V), q), \langle\, , \, \rangle_{4, k}\Big)$ is denoted by 
$g\ell _{q, W}^{\langle\, , \, \rangle_{4, k}}(V)$ and is called the {\bf $4$-th general linear Leibniz algebra} induced by $(q, W)$.

\medskip
We finish this section with the following 

\begin{definition}\label{def1.4} Let $W$ be a subspace of a vector space $V$ over a field $\mathbf{k}$ and let $q$ be a $W$-idempotent. 
\begin{description}
\item[(i)] A Lie algebra homomorphism $\varphi$ from a Lie algebra $(\mathcal{L},\,[ \,, \,] )$ to the $6$-th general linear Lie algebra 
$g\ell _{q, W}^{[ \,, \, ] _{6,k}}(V)$ is called a {\bf representation$^{6-th}$} of $\mathcal{L}$ on $V$ induced by $(q, W)$.
\item[(ii)] A Leibniz algebra homomorphism $\varphi$ from a Leibniz algebra 
$(\mathcal{L},\,\langle\, , \, \rangle )$ to the $4$-th general linear Leibniz algebra 
$g\ell _{q, W}^{\langle\, , \, \rangle_{4, k}}(V)$ is called a {\bf representation$^{4-th}$} of $\mathcal{L}$ on $V$ induced by $(q, W)$.
\end{description}
\end{definition}

\medskip
\section{Free Invariant Algebras}

\medskip
Let $(\hat{A}, \hat{q})$ be an invariant algebra, and let $\hat{i}$ be a map from a set $X$ to $(\hat{A}, \hat{q})$. The pair $\Big( (\hat{A}, \hat{q}), \hat{i}\Big)$ is called the {\bf free invariant algebra} generated by the set $X$ if the following universal property holds: given any invariant algebra $(A, q)$ and any map $\theta: X\to (A, q)$ there exists a unique invariant homomorphism $\hat{\theta}: (\hat{A}, \hat{q})\to (A, q)$ such that $\hat{\theta} \hat{i}=\theta$; that is, the following digram is commutative
$$
\begin{array}{ccc}
(\hat{A}, \hat{q})&\stackrel{\hat{\theta}}{\longrightarrow}&(A, q)\\
\shortstack{$\hat{i}$}\Bigg\uparrow&&\Bigg\uparrow\shortstack{$\theta$}\\&&\\
X&\stackrel{}{=}&X
\end{array}
$$

The free invariant algebra generated by a set is clearly unique. 

\medskip
We now construct free invariant algebras over a field $\mathbf{k}$. Let 
$X:=\{\, x_j \,|\, j\in J \,\}$ be a set, and let $\tilde{q}$ be a symbol which is not an element of $X$. Let $T(V)$ be the tensor algebra based on a vector space $V$, where 
$V=\displaystyle\bigoplus_{_{j\in J}}\mathbf{k}x_j \bigoplus \mathbf{k}\tilde{q}$ is the vector space over $\mathbf{k}$ with a basis $X\bigcup\{\tilde{q}\}$. Let $I$ be the ideal of $T(V)$ generated by
$$
\{\, \tilde{q}\otimes \tilde{q}-\tilde{q}, \,\, \tilde{q}\otimes a\otimes\tilde{q}-\tilde{q}\otimes a \,|\, a\in T(V) \,\}.
$$
Let $\hat{A}:=\displaystyle\frac{T(V)}{I}$ and $\hat{q}:=\tilde{q}+I$. Then $(\hat{A}, \hat{q})$ is an invariant algebra.

\medskip
\begin{proposition}\label{pr2.1} The pair $\Big( (\hat{A}, \hat{q}), \hat{i}\Big)$ is the free invariant algebra generated by the set $X=\{\, x_j \,|\, j\in J \,\}$, where the map
$\hat{i}: X\to \Big( (\hat{A}, \hat{q}), \hat{i}\Big)$ is defined by
$$
\hat{i}(x_j):=x_j+I \qquad\mbox{for $j\in J$.}
$$
\end{proposition}

\bigskip
Let $\hat{x}_j:=\hat{i}(x_j)=x_j+I$ for $j\in J$. The product of two elements $a$ and $b$ of 
$(\hat{A}, \hat{q})$ is simply denoted by $ab$. The next proposition gives a basic property of the free invariant algebra $\Big( (\hat{A}, \hat{q}), i\Big)$ generated by the set $X$.

\begin{proposition}\label{pr2.2} The following subset of  $\Big( (\hat{A}, \hat{q}), \hat{i}\Big)$ 
$$
\hat{S}:=\left\{\, \hat{1},\,\, \hat{q},\,\, \hat{x}_{j_1}\cdots\hat{x}_{j_m},\,\,
\hat{x}_{j_1}\cdots\hat{x}_{j_t}\hat{q}\hat{x}_{j_{t+1}}\cdots\hat{x}_{j_m} \,\left|\, 
\begin{array}{c}\{\,x_{j_1}, \cdots , x_{j_m}\,\}\subseteq X ,\\ m\in \mathcal{Z}_{\ge 1}, \,\, m\geq t\geq 0\, \end{array}\right.\,\right\}
$$
is a $\mathbf{k}$-basis of the vector space $(\hat{A}, \hat{q})$.
\end{proposition}

\section{Enveloping$^{6-th}$ Algebras of Lie Algebras}

\medskip
In the remaining part of this paper, we assume that {\bf the scalar $k$ is non-zero}.

\bigskip
Using invariant algebras and invariant homomorphisms in section 1, we introduce 
enveloping$^{6-th}$ algebras of Lie algebras in the following

\medskip
\begin{definition} Let $(\mathcal{L}, [\,, \,])$ be a Lie algebra (arbitrary dimensionality and characteristic). By a {\bf enveloping$^{6-th}$ algebra} of 
$(\mathcal{L}, [\,, \,])$ we will understand a pair $\Big((\mathcal{U}, \bar{q}), i\Big)$ composed of an invariant algebra $(\mathcal{U}, \bar{q})$ together with a map 
$i: \mathcal{L}\to (\mathcal{U}, \bar{q})$ satisfying the following two conditions:
\begin{description}
\item[(i)]  the map 
$i: (\mathcal{L}, [\,, \,])\to 
Lie\Big((\mathcal{U}, \bar{q}), [\,, \,]_{6, k}\Big)$ is a Lie algebra homomorphism; that is, $i$ is linear and
\begin{eqnarray*}
&&i([ x, y ])=[ i(x), i(y) ]_{6, k}=i(x)i(y)-i(y)i(x)+\\
&&\quad -i(x)i(y)\bar{q}+i(y)i(x)\bar{q}+ki(x)\bar{q}i(y)-ki(y)\bar{q}i(x)
\quad\mbox{for $x$, $y\in\mathcal{L}$,}
\end{eqnarray*}
\item[(ii)] given any invariant algebra $(A, q)$ and any Lie algebra homomorphism 
$f: (\mathcal{L}, [\,, \,])\to Lie\Big((A, q), [\,, \,]_{6,k}\Big)$ there
exists a unique invariant homomorphism 
$f': (\mathcal{U}, \bar{q})\to (A, q)$ such that $f=f'i$; that is, the following diagram is commutative:
$$
\begin{array}{ccc}
(\mathcal{U}, \bar{q})&\stackrel{f'}{\longrightarrow}&(A, q)\\
\shortstack{i}\Bigg\uparrow&&\Bigg\uparrow\shortstack{f}\\&&\\
\mathcal{L}&\stackrel{}{=}&\mathcal{L}.
\end{array}
$$
\end{description}
\end{definition}

Clearly, the enveloping$^{6-th}$ algebra of a Lie algebra $\mathcal{L}$, which is also denoted by 
$(\mathcal{U}^{6-th}_k(\mathcal{L}), \bar{q})$, is unique up to an invariant isomorphism.

\bigskip
We now construct the enveloping$^{6-th}$ algebra of a Lie algebra 
$(\mathcal{L}, [ \, ,\, ])$ over a field $\mathbf{k}$. 
Let $X=\{\, x_j \,|\, j\in J \,\}$ be a basis of 
$(\mathcal{L}, [\,, \,])$. Let $\Big( (\hat{A}, \hat{q}), \hat{i}\Big)$ be the free invariant algebra generated by the set $X$. By Proposition~\ref{pr2.2},
$\hat{X}:=\{\, \hat{x_j}:=\hat{i}(x_j) \,|\, j\in J \,\}$ is a linearly independent subset of $(\hat{A}, \hat{q})$. Hence, $\hat{i}$ can be extended to an injective linear map 
$\hat{i}: \mathcal{L}\to (\hat{A}, \hat{q})$. Let $\hat{x}:=\hat{i}(x)$ for all 
$x\in \mathcal{L}$. Let $R$ be the ideal of $(\hat{A}, \hat{q})$ which is generated by all the elements of the form
$$
\widehat{[ x, y ]}-\hat{x}\hat{y}+\hat{y}\hat{x}+\hat{x}\hat{y}\hat{q}
-\hat{y}\hat{x}\hat{q}-k\hat{x}\hat{q}\hat{y}+k\hat{y}\hat{q}\hat{x},
$$
where $\hat{x}$, $\hat{y}\in\mathcal{L}:=\hat{i}(\mathcal{L})$. Let 
\begin{equation}\label{eq2.10}
\mathcal{U}:=\displaystyle\frac{\hat{A}}{R}, \quad \bar{q}:=\hat{q}+R.
\end{equation}
Define a map $i: \mathcal{L}\to (\mathcal{U}, \bar{q})$ by
\begin{equation}\label{eq2.11}
i(x): =\hat{x}+R \qquad\mbox{for $x\in \mathcal{L}$}.
\end{equation}

We have the following

\begin{proposition}\label{pr2.3} The pair $\Big((\mathcal{U}, \bar{q}), i\Big)$ defined by (\ref{eq2.10}) and (\ref{eq2.11}) is an 
\linebreak enveloping$^{6-th}$ algebra for the Lie algebra 
$(\mathcal{L}, [\,, \,])$
\end{proposition}

\medskip
The following proposition gives the counterpart of P-B-W Theorem in the context of representations$^{6-th}$ of Lie algebras.

\medskip
\begin{proposition}\label{pr3.3} {\bf (The Extended$^{6-th}$ P-B-W Theorem)} Let $\mathcal{L}$ be a  Lie algebra with a basis  
$X:=\{\, x_j \,|\, j\in J\,\}$. If the set $J$ of indices is ordered, then
the following set of cosets 
$$
\hat{T}:=\left\{\begin{array}{c}
\hat{q}\hat{x}_{j_1}\cdots\hat{x}_{j_m}+R,\\
\hat{x}_{i_1}\cdots \hat{x}_{i_t}\hat{x}_{j_0}\hat{q}\hat{x}_{j_1}\cdots\hat{x}_{j_m}+R, \\
\hat{x}_{j_1}\cdots\hat{x}_{j_m}+R \\\end{array}
 \,\left|\, 
\begin{array}{c}x_{i_1}, \cdots , x_{i_t}, x_{j_0}, x_{j_1}, \cdots , x_{j_m}\in X,\\
i_t\geq\cdots \geq i_1,\\ 
j_m\ge \cdots\ge j_1\geq j_0,\\
t, m\in\mathcal{Z}_{\ge 0}
\end{array}\right.\right\}.
$$
is a basis for the enveloping$^{6-th}$ algebra $(\mathcal{U}, \bar{q}))$ of the Lie algebra $\mathcal{L}$.
\end{proposition}

\medskip
\section{Hopf-like$^{6-th}$ Algebras}

\medskip
Let $C$ be a vector space over a field $\mathbf{k}$. The canonical map $c\mapsto 1\otimes c$ is denoted by  $C\rightarrow\mathbf{k}\otimes C$, and the canonical map 
$c\mapsto c\otimes 1$ is denoted by $C\rightarrow C\otimes \mathbf{k}$  respectively, where $c\in C$. If $V$ and $W$ are vector spaces over a field $\mathbf{k}$, then the {\bf twist map} 
$$\tau : V\otimes W\to W\otimes V$$ 
is defined by $\tau (v\otimes w):=w\otimes v$, where $v\in V$ and $w\in W$.

\medskip
First, we introduce  the concept of a bialgebra with $\sigma$-counit.  

\begin{definition}\label{def5.3} An associative algebra $H$ over a field $\mathbf{k}$ is called a 
{\bf bialgebra with $\sigma$-counit} if there exist five linear maps $m: H\otimes H\to H$, 
$u: \mathbf{k}\to H$, $\Delta : H\to H\otimes H$,
$\varepsilon : H\to \mathbf{k}$ and $\sigma : H\to H$  such that the following eight diagrams are commutative:
\begin{itemize} 
\item associativity
$$
\begin{array}{ccc}
H\otimes H\otimes H
&\stackrel{m\otimes id}{\longrightarrow}&H\otimes H\\
id\otimes m\Bigg\downarrow&&\Bigg\downarrow m\\
H\otimes H&\stackrel{m}{\longrightarrow}&H
\end{array}
$$
\item unit
$$
\begin{array}{ccccc}
\mathbf{k}\otimes H
&\stackrel{u\otimes id}{\longrightarrow}&H\otimes H&
\stackrel{id\otimes u}{\longleftarrow}&H\otimes \mathbf{k}\\
\Bigg\downarrow&&\Bigg\downarrow m&&\Bigg\downarrow\\
H&\stackrel{id}{=}&H&\stackrel{id}{=}&H
\end{array}
$$
\item coassociativity
$$
\begin{array}{ccc}
H&\stackrel{\Delta}{\longrightarrow}&H\otimes H\\
\Delta\Bigg\downarrow&&\Bigg\downarrow id\otimes\Delta\\
H\otimes H&\stackrel{\Delta\otimes id}{\longrightarrow}&H\otimes H\otimes H
\end{array}
$$
\item $\sigma$-counit
$$
\begin{array}{ccccc}
\mathbf{k}\otimes H&\longleftarrow&H&\longrightarrow&H\otimes \mathbf{k}\\
\varepsilon\otimes id\Bigg\uparrow&&\Bigg\uparrow\sigma&&\Bigg\uparrow id\otimes \varepsilon\\
H\otimes H&\stackrel{\Delta}{\longleftarrow}&H&\stackrel{\Delta}{\longrightarrow}
&H\otimes H
\end{array}
$$
\item $\Delta$ is an algebra homomorphism
$$
\begin{array}{ccc}
H\otimes H&\stackrel{m}{\longrightarrow}H\stackrel{\Delta}{\longrightarrow}&H\otimes H\\
\Delta\otimes \Delta\Bigg\downarrow&&\Bigg\uparrow m\otimes m\\
H\otimes H\otimes H\otimes H&\stackrel{id\otimes \tau\otimes id}{\longrightarrow}&H\otimes H\otimes H\otimes H
\end{array}
$$
$$
\begin{array}{ccc}
H&\stackrel{\Delta}{\longrightarrow}&H\otimes H\\
u\Bigg\uparrow&&\Bigg\uparrow u\otimes u\\
\mathbf{k}&\longrightarrow&\mathbf{k}\otimes\mathbf{k}
\end{array}
$$
\item $\varepsilon$ is an algebra homomorphism
$$
\begin{array}{cc}
\begin{array}{ccc}
H\otimes H&\stackrel{\varepsilon\otimes\varepsilon}{\longrightarrow}&\mathbf{k}\otimes\mathbf{k}\\
m\Bigg\downarrow&&\Bigg\downarrow \\
H&\stackrel{\varepsilon}{\longrightarrow}&\mathbf{k}
\end{array}\qquad\qquad&
\begin{array}{ccc}
H&=&H\\
u\Bigg\uparrow&&\Bigg\downarrow \varepsilon\\
\mathbf{k}&=&\mathbf{k}
\end{array}
\end{array}
$$
\end{itemize}
The maps $\Delta$ and $\varepsilon$ are called the 
{\bf comultiplication} and the {\bf $\sigma$-counit} respectively. A bialgebra H with 
$\sigma$-counit is also denoted by $H{m, u, \Delta \choose  \varepsilon , \sigma }$ if the five linear maps have to be indicated explicitly.
\end{definition}

Based on bialgebras with $\sigma$-counit, we now introduce Hopf-like algebras in the following 

\begin{definition}\label{def4.2} A  bialgebra $H=H{m, u, \Delta \choose  \varepsilon , \sigma }$  with $\sigma$-counit over a field $\mathbf{k}$ is called a 
{\bf Hopf-like$^{6-th}$ algebra} if there exists a linear map $S : H\to H$ such that the following diagram is commutative:
\begin{equation}\label{eq588}
\begin{array}{ccccccc}
H&\stackrel{\Delta}{\longrightarrow}&H\otimes H&\stackrel{S\otimes id}{\longrightarrow}
&H\otimes H&\stackrel{m}{\longrightarrow}&H\\
\Big|\Big|&&&&&&\Big|\Big|\\
H&\stackrel{S}{\longrightarrow}&H&\stackrel{\varepsilon}{\longrightarrow}&\mathbf{k}
&\stackrel{u}{\longrightarrow}&H\\
\Big|\Big|&&&&&&\Big|\Big|\\
H&\stackrel{\Delta}{\longrightarrow}&H\otimes H&\stackrel{id\otimes S}{\longrightarrow}
&H\otimes H&\stackrel{m}{\longrightarrow}&H
\end{array}
\end{equation}
\end{definition}
The linear map $S$ satisfying (\ref{eq588}) is called the {\bf antipode$^{6-th}$-like} of the Hopf-like$^{6-th}$ algebra.

\medskip
The first generalization of the Hopf algebra structure for the ordinary enveloping algebra of a Lie algebra is given in the following

\begin{proposition}\label{pr5.4} The enveloping$^{6-th}$ algebra $U^{6-th}_k(\mathcal{L})$ of a Lie algebra $\mathcal{L}$ is a Hopf-like$^{6-th}$ algebra, where the comultiplication $\Delta$, the $\sigma$-counit 
$\varepsilon$ and the antipode$^{6-th}$-like $S$ satisfy
$$
\Delta (q)=q\otimes q,
$$
\begin{eqnarray*}
\Delta (x)&=&(x+kqx-xq)\otimes 1+1\otimes (x+kqx-xq)+\nonumber\\
&& +(1-k)qx\otimes q+(1-k)q\otimes qx \qquad\mbox{for $x\in \mathcal{L}$,}
\end{eqnarray*}
$$
\varepsilon (q)=1,\qquad \varepsilon (x)=0 \qquad\mbox{for $x\in \mathcal{L}$,}
$$
$$
S(q)=1-q, \qquad 
S(x)=-\frac1k x-kqx+\frac1k xq\qquad\mbox{for $x\in \mathcal{L}$}
$$
and the linear map $\sigma : U^{6-th}_k(\mathcal{L})\to U^{6-th}_k(\mathcal{L})$ is defined by 
$$\sigma (a):=a+qa-aq\qquad\mbox{for $a\in \mathcal{U}^{6-th}_k(\mathcal{L})$.}$$
\end{proposition}

\medskip
\section{Enveloping$^{\langle 4-th \rangle}$ algebras of Leibniz algebras}

The concept of the enveloping$^{\langle 4-th \rangle}$ algebra of a Leibniz Algebra is introduced in the following

\begin{definition}\label{def3.1} Let $(\mathcal{L}, \langle, \rangle)$ be a Leibniz algebra (arbitrary dimensionality and characteristic). By a {\bf enveloping$^{\langle 4-th \rangle}$ algebra} of 
$(\mathcal{L}, \langle, \rangle)$ we will understand a pair $\Big((\mathcal{U}, \bar{q}), i\Big)$ composed of an invariant algebra $(\mathcal{U}, \bar{q})$ together with a map 
$i: \mathcal{L}\to (\mathcal{U}, \bar{q})$ satisfying the following two conditions:
\begin{description}
\item[(i)]  the map 
$i: (\mathcal{L}, \langle, \rangle)\to Leib\Big((\mathcal{U}, \bar{q}), \langle, \rangle_{4,k}\Big)$ is a Leibniz algebra homomorphism; that is, $i$ is linear and
\begin{eqnarray*}
&&i(\langle x, y\rangle)=\langle i(x), i(y)\rangle_{4,k}=i(x)i(y)-i(y)i(x)+\\
&&\quad +i(y)\bar{q}i(x)-i(x)i(y)\bar{q}+ki(x)\bar{q}i(y)-k\bar{q}i(x)i(y)
\quad\mbox{for $x$, $y\in\mathcal{L}$,}
\end{eqnarray*}
\item[(ii)] given any invariant algebra $(A, q)$ and any Leibniz algebra homomorphism 
$f: (\mathcal{L}, \langle, \rangle)\to Leib\Big((A, q), \langle, \rangle_{4,k}\Big)$ there
exists a unique invariant homomorphism $f': (\mathcal{U}, \bar{q})\to (A, q)$ such that $f=f'i$; that is, the following diagram is commutative:
$$
\begin{array}{ccc}
(\mathcal{U}, \bar{q})&\stackrel{f'}{\longrightarrow}&(A, q)\\
\shortstack{i}\Bigg\uparrow&&\Bigg\uparrow\shortstack{f}\\&&\\
\mathcal{L}&\stackrel{}{=}&\mathcal{L}.
\end{array}
$$
\end{description}
\end{definition}

\medskip
Clearly, the enveloping$^{\langle 4-th \rangle}$ algebra of a Leibniz algebra $\mathcal{L}$, which is also denoted by $(\mathcal{U}^{\langle 4-th \rangle}_k(\mathcal{L}), \bar{q})$, is unique up to an invariant isomorphism.

\medskip
We now construct the enveloping invariant algebra of a Leibniz algebra 
$(\mathcal{L}, \langle, \rangle)$ over a field $\mathbf{k}$. Let $X=\{\, x_j \,|\, j\in J \,\}$ be a basis of $(\mathcal{L}, \langle, \rangle)$. Let $\Big( (\hat{A}, \hat{q}), \hat{i}\Big)$ be the free invariant algebra generated by the set $X$. By Proposition~\ref{pr2.2},
$\hat{X}:=\{\, \hat{x}_j:=\hat{i}(x_j) \,|\, j\in J \,\}$ is a linearly independent subset of $(\hat{A}, \hat{q})$. Hence, $i$ can be extended to an injective linear map 
$\hat{i}: \mathcal{L}\to (\hat{A}, \hat{q})$. Let $\hat{x}:=\hat{i}(x)$ for all 
$x\in \mathcal{L}$. Let $R$ be the ideal of $(\hat{A}, \hat{q})$ which is generated by all the elements of the form
$$
\widehat{\langle x, y\rangle}-\hat{x}\hat{y}+\hat{y}\hat{x}-\hat{y}\hat{q}\hat{x}+\hat{x}\hat{y}\hat{q}-k\hat{x}\hat{q}\hat{y}+k\hat{q}\hat{y}\hat{x},
$$
where $x$, $y\in\mathcal{L}$ and we identify $\mathcal{L}$ with $\hat{i}(\mathcal{L})$. Let 
\begin{equation}\label{eq2.8}
\mathcal{U}:=\displaystyle\frac{\hat{A}}{R}, \quad \bar{q}:=\hat{q}+R.
\end{equation}
Define a map $i: \mathcal{L}\to (\mathcal{U}, \bar{q})$ by
\begin{equation}\label{eq2.9}
i(x): =\hat{x}+R \qquad\mbox{for $x\in \mathcal{L}$}.
\end{equation}

\medskip
\begin{proposition}\label{pr2.3} The pair $\Big((\mathcal{U}, \bar{q}), i\Big)$ defined by (\ref{eq2.8}) and (\ref{eq2.9}) is the enveloping$^{\langle 4-th \rangle}$ algebra of the Leibniz algebra 
$(\mathcal{L}, \langle, \rangle)$
\end{proposition}

\medskip
The following proposition gives the counterpart of P-B-W Theorem in the context of representations$^{\langle 4-th \rangle}$ of Leibniz algebras.

\begin{proposition}{\bf (The Extended$^{\langle 4-th \rangle}$ P-B-W Theorem)} Let $\mathcal{L}$ be a Leibniz algebra with a basis  
$X:=\{\, x_j \,|\, j\in J\,\}$. If the set $J$ of indices is ordered, 
then the following set of cosets of model monomials
$$
T:=\left\{\,\begin{array}{c}
\hat{q}\hat{x}_{j_1}\cdots\hat{x}_{j_m}+R, \\
\hat{x}_j\hat{q}\hat{x}_{j_1}\cdots\hat{x}_{j_m}+R, \\
\hat{x}_{j_1}\cdots\hat{x}_{j_m}+R \\\end{array}
 \,\left|\, 
\begin{array}{c}\{\,x_{j}, x_{j_1}, \cdots , x_{j_m}\,\}\subseteq X ,\\
\mbox{$j_m\ge \cdots \ge j_1$ and $m\in\mathcal{Z}_{\ge 0}$}\end{array}\right.\,\right\}
$$
is a basis for the enveloping$^{\langle 4-th \rangle}$ algebra of the Leibniz algebra $\mathcal{L}$.
\end{proposition}

\medskip
\section{Hopf-like$^{\langle 4-th \rangle}$  Algebras}

We begin this section by introducing two generalizations of Hopf algebras.

\begin{definition}\label{def5.3} Let $H=(H, m, u, \Delta ,  \varepsilon )$  be a bialgebra, where  $(H, m, u)$ gives the associative algebra structure of $H$, and 
$(H, \Delta ,  \varepsilon )$ gives the coalgebra structure of $H$.
\begin{description}
\item[(i)] H is called a {\bf Hopf-like$^{{\langle 4-th \rangle}_1}$  algebra} if there exists a linear map $\mathcal{S} : H\to H$ such that the following diagram is commutative:
\begin{equation}\label{eq58}
\begin{array}{ccccccc}
H&\stackrel{\Delta}{\longrightarrow}&H\otimes H&\stackrel{\mathcal{S}\otimes id}{\longrightarrow}
&H\otimes H&\stackrel{m}{\longrightarrow}&H\\
\Big|\Big|&&&&&&\Big|\Big|\\
H&\stackrel{\mathcal{S}}{\longrightarrow}&H&\stackrel{\varepsilon}{\longrightarrow}&\mathbf{k}
&\stackrel{u}{\longrightarrow}&H\\
\Big|\Big|&&&&&&\Big|\Big|\\
H&\stackrel{\Delta}{\longrightarrow}&H\otimes H&\stackrel{id\otimes \mathcal{S}}{\longrightarrow}
&H\otimes H&\stackrel{m}{\longrightarrow}&H
\end{array}
\end{equation}
The linear map $\mathcal{S}$ satisfying  (\ref{eq58}) is called the 
{\bf antipode$^{{\langle 4-th \rangle}_1}$-like} of the Hopf-like$^{{\langle 4-th \rangle}_1}$ algebra $H$.
\item[(ii)] H is called a {\bf Hopf-like$^{{\langle 4-th \rangle}_2}$  algebra} if there exist an algebra homomorphism $\sigma : H\to H$ and a linear map $\mathcal{S} : H\to H$ such that the following diagram is commutative:
\begin{equation}\label{eq59}
\begin{array}{ccccccc}
H&\stackrel{\Delta}{\longrightarrow}&H\otimes H&\stackrel{\mathcal{S}\otimes id}{\longrightarrow}
&H\otimes H&\stackrel{m}{\longrightarrow}H\stackrel{\sigma}{\longrightarrow}&H\\
\Big|\Big|&&&&&&\Big|\Big|\\
H&\stackrel{\mathcal{S}}{\longrightarrow}&H&\stackrel{\varepsilon}{\longrightarrow}&\mathbf{k}
&\stackrel{u}{\longrightarrow}&H\\
\Big|\Big|&&&&&&\Big|\Big|\\
H&\stackrel{\Delta}{\longrightarrow}&H\otimes H&\stackrel{id\otimes \mathcal{S}}{\longrightarrow}
&H\otimes H&\stackrel{m}{\longrightarrow}&H
\end{array}
\end{equation}
The linear map $\mathcal{S}$ satisfying (\ref{eq59}) is called the 
{\bf antipode$^{{\langle 4-th \rangle}_2}$-like} of the Hopf-like$^{{\langle 4-th \rangle}_2}$ algebra $H$.
\end{description}
\end{definition}

\medskip
Clearly, a Hopf algebra is both a Hopf-like$^{{\langle 4-th \rangle}_1}$ algebra and a 
Hopf-like$^{{\langle 4-th \rangle}_2}$ algebra with $\sigma =id$. 

\medskip
The other two generalizations of the Hopf algebra structure for the ordinary enveloping algebra of a Lie algebra is given in the following

\begin{proposition} Let $(\mathcal{U}^{\langle 4-th \rangle}_k(\mathcal{L}), q)$ be the enveloping$^{\langle 4-th \rangle}$ algebra of a Leibniz algebra $\mathcal{L}$.
\begin{description}
\item[(i)] If $k=1$, then $(\mathcal{U}^{\langle 4-th \rangle}_1(\mathcal{L}), q)$ is a 
Hopf-like$^{{\langle 4-th \rangle}_1}$ algebra;
\item[(ii)] If $k\ne 0$, then $(\mathcal{U}^{\langle 4-th \rangle}_k(\mathcal{L}), q)$ is a 
Hopf-like$^{{\langle 4-th \rangle}_2}$ algebra,
\end{description}
where the comultiplication $\Delta$,
the counit $\varepsilon$, the antipode-like$^{{\langle 4-th \rangle}_i}$ $\mathcal{S}$ with $i=1$, $2$ satisfy
$$
\Delta (q)=q\otimes q,
$$
\begin{eqnarray*}
\Delta (x)&=&(x+kqx-xq)\otimes 1+1\otimes (x+kqx-xq)+\nonumber\\
&& +(-kqx+xq)\otimes q+q\otimes (-kqx+xq) \qquad\mbox{for $x\in \mathcal{L}$,}
\end{eqnarray*}
$$
\varepsilon (q)=1, \qquad  \varepsilon (x)=0 \qquad\mbox{for $x\in \mathcal{L}$,}
$$
$$
\mathcal{S}(q)=1-q,\qquad 
\mathcal{S}(x)=-\frac1k x+\Big(\frac1k -k\Big) xq\qquad\mbox{for $x\in \mathcal{L}$.}
$$
and the algebra homomorphism 
$\sigma $ is defined by 
$$\sigma (a):=a+qa-aq\qquad\mbox{for $a\in \mathcal{U}^{\langle 4-th \rangle}_k(\mathcal{L})$.}$$
\end{proposition}

\end{document}